\documentclass[10pt]{article}
\usepackage{amsfonts}
\usepackage{latexsym}
\usepackage{cite}
\usepackage{amsmath,amsfonts,latexsym,amssymb}
\usepackage[mathscr]{eucal}
\usepackage{cases,color}
\usepackage{amsthm}

\usepackage[bf,small]{caption2}
\usepackage{float}
\usepackage{graphicx}
\usepackage{amsmath}
\usepackage{amssymb}
\usepackage[all]{xy}

\newtheorem{theorem}{theorem}[section]

\newtheorem{conv}[theorem]{Convention}

\newtheorem{exmp}[theorem]{Example}
\newtheorem{lem}[theorem]{Lemma}

\newtheorem{prop}[theorem]{Proposition}

\newtheorem{rmk}[theorem]{Remark}
\newtheorem{thm}[theorem]{Theorem}

\begin{document}

\title{\textbf{The ${\rm SL}(2,\mathbb{C})$-character variety of a Montesinos knot}}
\author{\Large Haimiao Chen}

\date{}
\maketitle

\begin{abstract}
  For each Montesinos knot $K$, we propose an efficient method to explicitly determine the irreducible ${\rm SL}(2,\mathbb{C})$-character variety, and show that it can be decomposed as $\mathcal{X}_0(K)\sqcup\mathcal{X}_1(K)\sqcup\mathcal{X}_2(K)\sqcup\mathcal{X}'(K)$, where $\mathcal{X}_0(K)$ consists of trace-free characters, $\mathcal{X}_1(K)$ consists of characters of ``unions" of representations of rational knots (or rational link, which appears at most once), $\mathcal{X}_2(K)$ is an algebraic curve, and $\mathcal{X}'(K)$ consists of finitely many points when $K$ satisfies a generic condition.

\end{abstract}

\section{Introduction}

Let $G={\rm SL}(2,\mathbb{C})$ and let $\Gamma$ be a finitely presented group. A {\it $G$-representation} of $\Gamma$ is a homomorphism $\rho:\Gamma\to G$.
The {\it $G$-representation variety} of $\Gamma$ is $\mathcal{R}_G(\Gamma):=\hom(\Gamma,G)$.
Call a representation $\rho$ {\it reducible} if elements in ${\rm Im}(\rho)$ have a common eigenvector, which is equivalent to ${\rm Im}(\rho)\subset \mathbf{a}U\mathbf{a}^{-1}$ for some $\mathbf{a}\in G$, where
$U$ denotes the subgroup of $G$ consisting of upper-triangular matrices;
in particular, call $\rho$ {\it abelian} if ${\rm Im}(\rho)$ is abelian. Call $\rho$ {\it irreducible} if $\rho$ is not reducible.
Let $\mathcal{R}^{\rm irr}_G(\Gamma)$ denote the space of irreducible representations.

The {\it character} of $\rho$ is by definition the function $\chi_\rho:\Gamma\to\mathbb{C}$, $g\mapsto{\rm tr}(\rho(g)).$
As a fact (\cite{CS83} Proposition 1.5.2), two irreducible representations $\rho,\rho'$ are conjugate (meaning that there exists $\mathbf{a}\in G$ such that $\rho'(x)=\mathbf{a}\rho(x)\mathbf{a}^{-1}$ for all $x\in\Gamma$) if and only if $\chi_\rho=\chi_{\rho'}$.
Focusing on irreducible representations, the (irreducible) {\it $G$-character variety} of $\Gamma$ is
$\mathcal{X}^{\rm irr}_G(\Gamma)=\{\chi_\rho\colon\rho\in\mathcal{R}^{\rm irr}_G(\Gamma)\}.$

When $\Gamma=\pi(L):=\pi_1(S^3-L)$ for a link $L$, we call $\mathcal{X}^{\rm irr}_G(\Gamma)$ the $G$-character variety of $L$, and so on.
Throughout the paper, we usually work directly on representations, which are more convenient. A representation whose conjugacy ambiguity is fixed will always be identified with its character.

There are many reasons for caring about representations. Here are two of them.
First, in general, having a representation $\rho:\pi(L)\to {\rm GL}(d,\mathbb{F})$, one can define the {\it twisted Alexander polynomial} of $L$ associated to $\rho$, which is useful in many situations (referred to \cite{FV08} and the references therein).
Second, nowadays a popular topic is the {\it left orderability} of $\pi_1(M)$ for 3-manifolds $M$; here a group $G$ is called left-orderable if it admits a total ordering $<$ such that $g<h$ implies $fg<fh$ for all $f,g,h\in G$. For a 3-manifold $M$ resulting from a Dehn surgery on a knot $K$, a systematic approach to understanding the left orderability of $\pi_1(M)$ was proposed by Culler-Dunfield \cite{CD18}, by using continuous paths of ${\rm SL}(2,\mathbb{R})$-representations of $K$.

Computing character varieties of links is a notoriously difficult problem. The importance for computation lies in that knowledge on character variety is still too limited, and concrete results may give directions for further research.
Till now, character varieties have been determined for only a few links. Existing results include:
torus knots \cite{Mu09}, double twist knots \cite{MPL11}, double twist links \cite{PT15}, $(-2,2m+1,2n)$-pretzel links and twisted Whitehead links \cite{Tr16}, and classical pretzel knots \cite{Ch18-2,Ch19}. Recently, the result of \cite{Ch18-2} was applied by Khan and Tran \cite{KT20} to show left orderability for some Dehn surgeries on odd pretzel knots.

As a related work, in \cite{PP13} Paoluzzi and Porti studied the ${\rm SL}(2,\mathbb{C})$-character varieties of Montesinos knots of {\it Kinoshita-Terasaka type} using the ``bending trick", to reveal some interesting phenomena; for instance, the character variety contains a $d$-dimensional component for each $d$ in a certain range.

In this article, we show
\begin{thm} \label{thm:main}
For each Montesinos knot $K$, the irreducible character variety has a decomposition
$\mathcal{X}^{\rm irr}_{{\rm SL}(2,\mathbb{C})}(K)=\mathcal{X}_0(K)\sqcup\mathcal{X}_1(K)\sqcup\mathcal{X}_2(K)\sqcup\mathcal{X}'(K),$
where
\begin{itemize}
  \item $\mathcal{X}_0(K)$ consists of trace-free characters;
  \item $\mathcal{X}_1(K)$ consists of characters of ``unions" of representations of rational knots or maybe one rational link;
  \item $\mathcal{X}_2(K)$ is a high-genus algebraic curve, whose defining equation can be written down;
  \item $\mathcal{X}'(K)$ is a finite set when $K$ satisfies a generic condition.
\end{itemize}
\end{thm}
The ``generic condition" for $K$ is precisely stated in Proposition \ref{prop:X'}.

Actually we give an efficient method for determining the character variety. The techniques used are elementary, but the results are amazing. Each of the four parts $\mathcal{X}_0(K),\mathcal{X}_1(K),\mathcal{X}_2(K),\mathcal{X}'(K)$ has distinguishable feature.
When the number $m$ of strands of $K$ is larger than 3, many new phenomena appear, in contrast to $m=3$. $\mathcal{X}_1(K)$ is concerned with representations from $F_m$ (the free group on $m$ generators) to ${\rm SL}(2,\mathbb{C})$. The last part $\mathcal{X}'(K)$ is rather subtle, and can be regarded as certain kind of singularity, which we manage to handle.

This is the first time to systematically deal with a class of knots whose fundamental groups can have arbitrarily many generators.
The structural description of $\mathcal{X}^{\rm irr}_{{\rm SL}(2,\mathbb{C})}(K)$ provides necessary conditions for deciding whether a given knot is Montesinos, which in general is a difficult problem \cite{IJ13}.
Explicit descriptions will provide plenty of examples on character varieties; in particular, examples of character varieties having high-dimensional components.

The paper is organized as follows. Section 2 is a collection of basic notions. In Section 3, representations of rational tangles are investigated in detail; some fine things are clarified and several useful formulas are derived. In Section 4, a method for finding all irreducible representations for each Montesinos knot is developed, and then Theorem \ref{thm:main} is established. We assume that the values of characters taken at meridians are nonzero, since the case when the values vanish had been worked out in \cite{Ch18-1}.

\section{Preparation}

\subsection{Tangles and links}

We adopt the notations and conventions used in \cite{Ch18-1}.

By a {\it tangle} we simultaneously mean an un-oriented tangle diagram and the 1-dimensional manifold embedded in $\mathbb{R}^3$ it stands for.

\begin{figure} [h]
  \centering
  \includegraphics[width=0.65\textwidth]{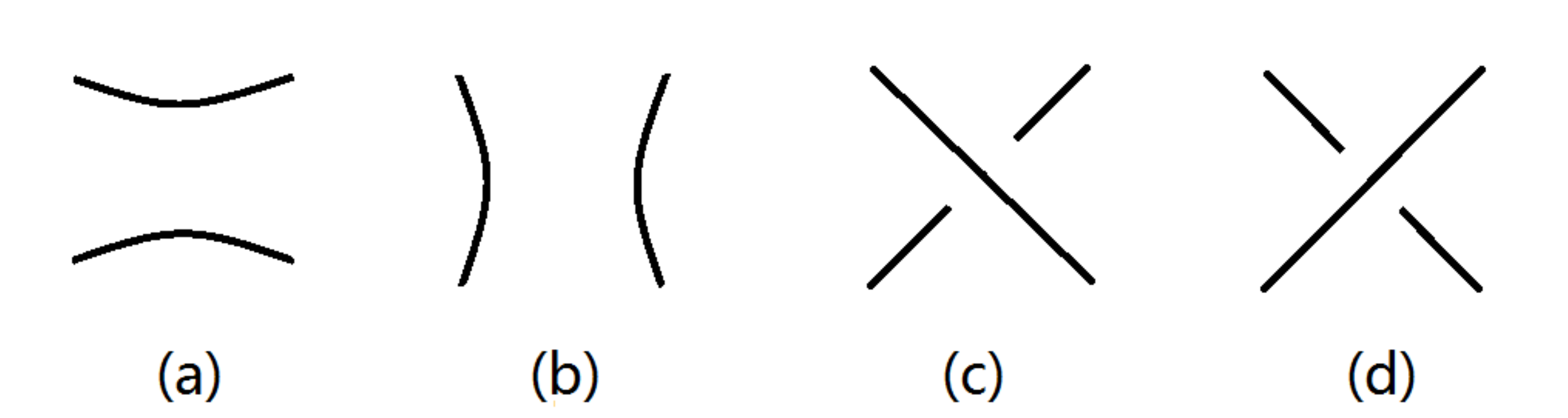}\\
  \caption{Four elements in $\mathcal{T}_2^2$: (a) $[0]$, (b) $[\infty]$, (c) $[1]$, (d) $[-1]$.}\label{fig:basic}
\end{figure}

\begin{figure}[h]
  \centering
  \includegraphics[width=0.55\textwidth]{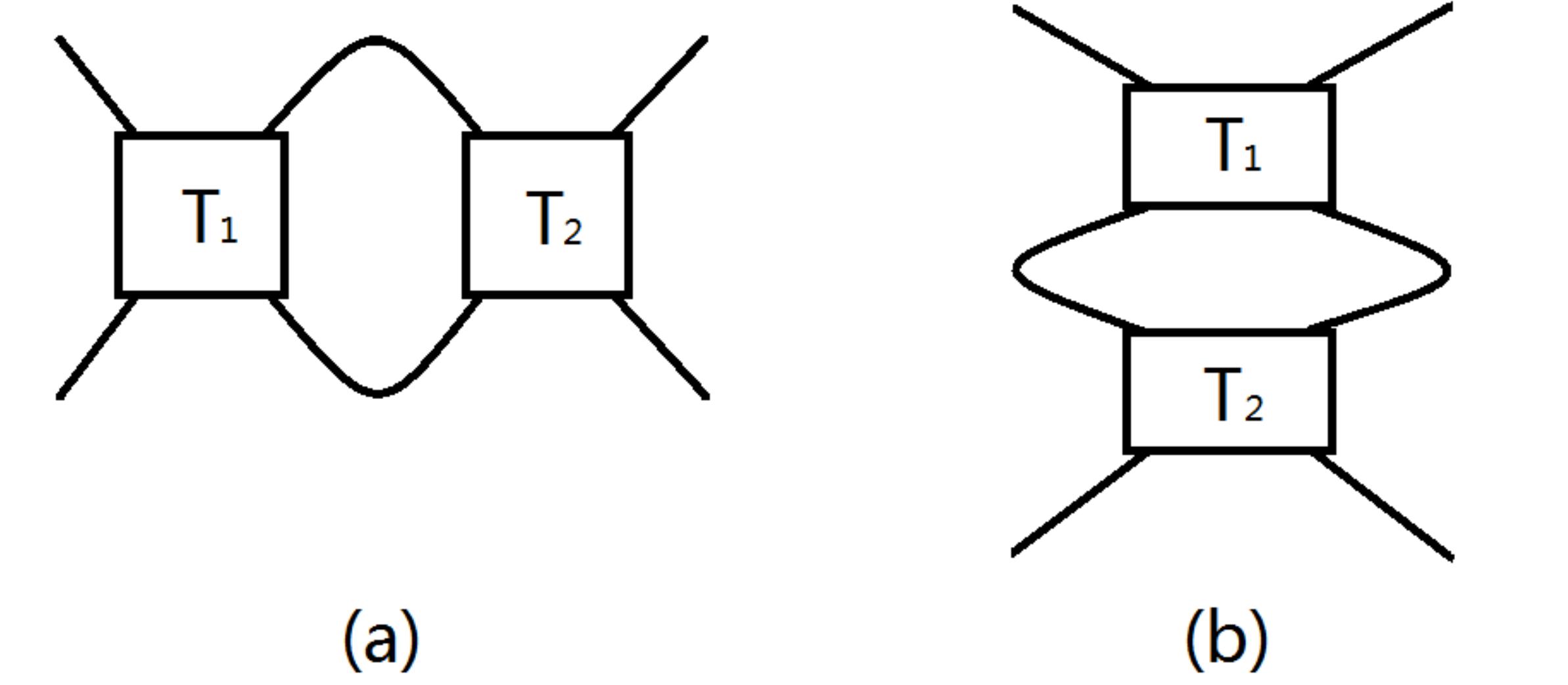}\\
  \caption{(a) $T_{1}+T_{2}$; (b) $T_{1}\ast T_{2}$.} \label{fig:composition}
\end{figure}

Let $\mathcal{T}_2^2$ denote the set of 4-ended tangles. 
Four of the simplest elements of $\mathcal{T}_2^2$ are $[0]$, $[\infty]$, $[1]$, $[-1]$, as shown in Figure \ref{fig:basic}. Defined on $\mathcal{T}_2^2$ are horizontal composition $+$ and vertical composition $\ast$, as illustrated in Figure \ref{fig:composition}.

For $k\ne 0$, the horizontal composite of $|k|$ copies of $[1]$ (resp. $[-1]$) is denoted by $[k]$ if $k>0$ (resp. $k<0$), and the
vertical composite of $|k|$ copies of $[1]$ (resp. $[-1]$) is denoted by $[1/k]$ if $k>0$ (resp. $k<0$).
Given $k_{1},\ldots,k_{s}\in\mathbb{Z}$, the {\it rational tangle} $[[k_{1}],\ldots,[k_{s}]]$ is defined as
$$\begin{cases}
[k_{1}]\ast [1/k_{2}]+\cdots+[k_{s}], &2\nmid s, \\
[1/k_{1}]+[k_{2}]\ast\cdots +[k_{s}], &2\mid s.
\end{cases}$$
The {\it fraction} of the tangle $[[k_{1}],\ldots,[k_{s}]]$ is given by the continued fraction $[[k_{1},\ldots,k_{s}]]\in\mathbb{Q}$, which is defined inductively as
$$[[k_{1}]]=k_{1}; \qquad
[[k_{1},\ldots,k_{j}]]=k_{j}+1/[[k_{1},\ldots,k_{j-1}]], \quad 2\le j\le s.$$
Denote the tangle $[[k_{1}],\ldots,[k_{s}]]$ as $[p/q]$ when its fraction equals $p/q$, where $p$ is coprime to $q$, and call it {\it odd} (resp. {\it even}) if $p$ is odd (resp. even).

\begin{figure}[h]
  \centering
  \includegraphics[width=6cm]{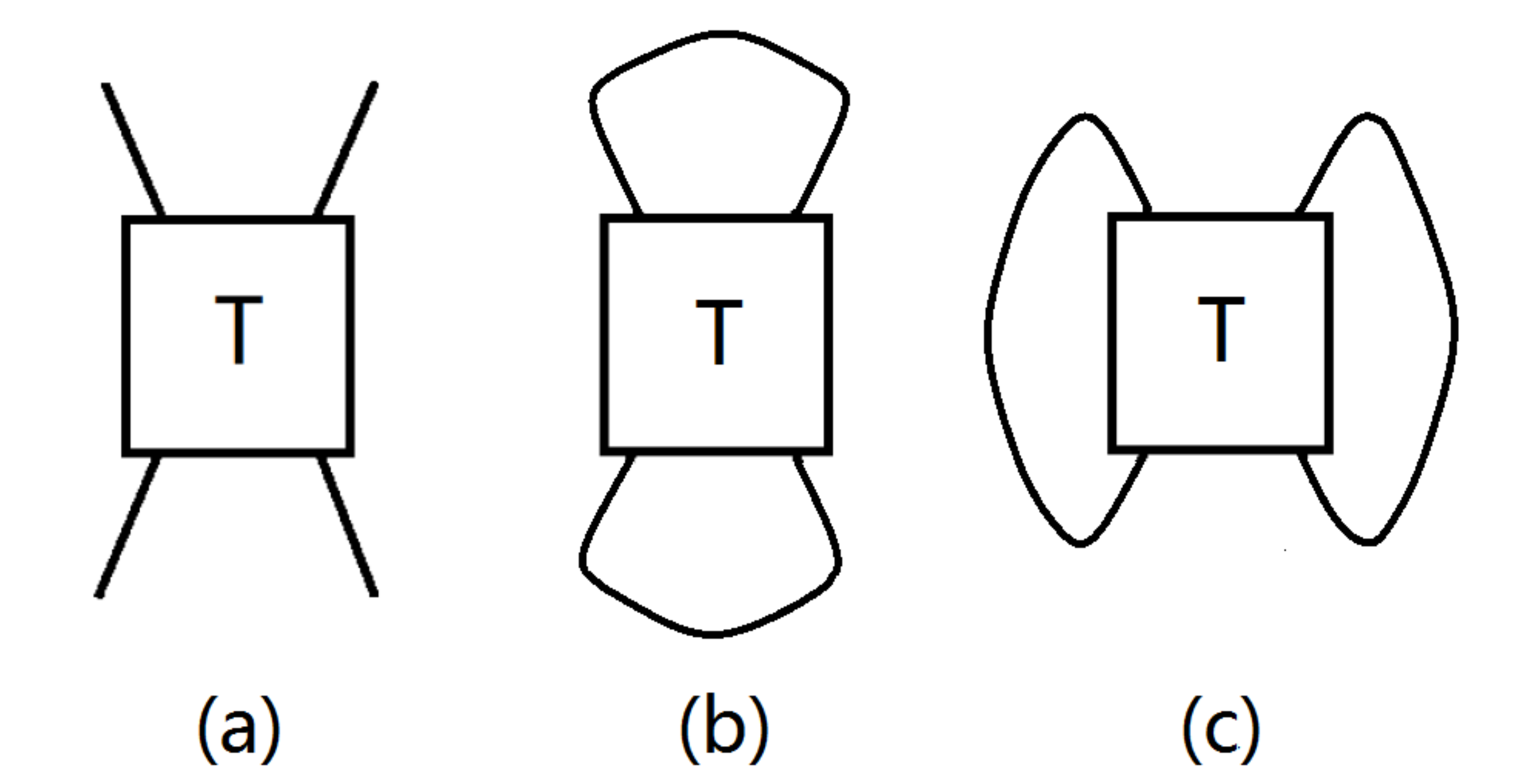}\\
  \caption{(a) a tangle $T$; (b) the numerator $N(T)$; (c) the denominator $D(T)$.} \label{fig:D-N}
\end{figure}

Each $T\in\mathcal{T}_2^2$ gives rise to two links: the {\it numerator} $N(T)$ and the {\it denominator} $D(T)$, as shown in Figure \ref{fig:D-N}. For $p/q\in\mathbb{Q}$, the link $N([p/q])$ is called a {\it rational link}
or {\it 2-bridge link}, which has one component (resp. two components) if $p$ is odd (resp. even). A link of the form $D([p_1/q_1]\ast\cdots\ast[p_m/q_m])$ is called a {\it Montesinos link} and denoted by $M(p_1/q_1,\ldots,p_m/q_m)$;
in particular, $M(p_1,\ldots,p_m)$ is a {\it pretzel link}.

\begin{conv}
\rm An Montesinos knot $M(p_1/q_1,\ldots,p_m/q_m)$ is called {\it odd} if $2\nmid p_i$ for all $i$, and called {\it even} if $2\mid p_i$ for exactly one $i$. In the second case, by cyclically permuting the indices, we always assume $2\mid p_m$.
Even Montesinos knots are the same as those of Kinoshita-Terasaka type. 

Note that the northwest and northeast end of $[p_i/q_i]$ belong to the same component if and only if $M(p_1/q_1,\ldots,p_m/q_m)$ is even and $i=m$.
\end{conv}

\subsection{On $2\times 2$ matrices}

We use bold letters to denote elements of ${\rm SL}(2,\mathbb{C})$. 

Introduce some special elements:
\begin{align*}
\mathbf{e}=\left(\begin{array}{cc} 1 & 0 \\ 0 & 1 \end{array}\right), \qquad
\mathbf{p}=\left(\begin{array}{cc} 1 & 1 \\ 0 & 1 \end{array}\right), \qquad
\mathbf{w}=\left(\begin{array}{cc} 0 & 1 \\ -1 & 0 \end{array}\right),  \\
\mathbf{d}(\kappa)=\left(\begin{array}{cc} \kappa & 0 \\ 0 & \kappa^{-1} \end{array}\right), \quad
\mathbf{u}^+_\kappa(\xi)=\left(\begin{array}{cc} \kappa & \xi \\ 0 & \kappa^{-1} \end{array}\right), \quad
\mathbf{u}_\kappa^-(\xi)=\left(\begin{array}{cc} \kappa & 0 \\ \xi & \kappa^{-1}\end{array}\right).
\end{align*}
The {\it Weyl element} $\mathbf{w}$ has the property that $\mathbf{w}\mathbf{x}\mathbf{w}^{-1}=(\mathbf{x}')^{-1}$ for each $\mathbf{x}$, where $\mathbf{x}'$ denotes the transpose of $\mathbf{x}$.

Let $M$ denote the set of $2\times 2$ matrices. Let $U$ denote the subgroup of ${\rm SL}(2,\mathbb{C})$ consisting of upper-triangular elements.

For $\lambda$ with $\lambda+\lambda^{-1}\ne t^2-2,\pm 2$ and $\mu\ne 0$, put
\begin{align*}
\mathbf{h}_{t}^{\lambda}(\mu)=
\frac{1}{\lambda+1}\left(\begin{array}{cc} \lambda t & \mu \\ (t^2-\lambda-\lambda^{-1}-2)\lambda\mu^{-1} & t \end{array}\right).
\end{align*}
For $t\ne 0$, put
\begin{align*}
\mathbf{k}_{t}(\alpha)=\left(\begin{array}{cc} t/2+\alpha & (t^2/4-1-\alpha^2)/(2t) \\ 2t & t/2-\alpha \end{array}\right).
\end{align*}

Direct computation leads to
\begin{align}
{\rm tr}\big(\mathbf{h}_{t}^{\lambda}(\mu)^{-1}\mathbf{h}_{t}^{\lambda}(\nu)\big)&=\frac{2t^2+(\lambda+\lambda^{-1}+2-t^2)(\mu\nu^{-1}+\mu^{-1}\nu)}{\lambda+\lambda^{-1}+2},  \label{eq:tr-h} \\
{\rm tr}\big(\mathbf{k}_t(\alpha)^{-1}\mathbf{k}_t(\beta)\big)&=2+(\alpha-\beta)^2. \label{eq:tr-k}
\end{align}

The reason for introducing these matrices is clear from the following lemma.
\begin{lem}\label{lem:key}
Suppose $\mathbf{a}_1,\mathbf{a}_2\in{\rm SL}(2,\mathbb{C})$ with ${\rm tr}(\mathbf{a}_1)={\rm tr}(\mathbf{a}_2)=t$.
\begin{enumerate}
  \item[\rm(a)] If $\mathbf{a}_1\mathbf{a}_2=\mathbf{d}(\lambda)$ with $\lambda+\lambda^{-1}\ne t^2-2, \pm 2$,
                then $\mathbf{a}_1=\mathbf{h}_{t}^{\lambda}(\mu)$ and $\mathbf{a}_2=\mathbf{h}_{t}^{\lambda}(-\lambda^{-1}\mu)$ for some $\mu\ne 0$.
  \item[\rm(b)] If $\mathbf{a}_1\mathbf{a}_2=\mathbf{p}$, then $\mathbf{a}_1, \mathbf{a}_2\in U$.
  \item[\rm(c)] If $\mathbf{a}_1\mathbf{a}_2=-\mathbf{p}$ and $t\ne 0$, then $\mathbf{a}_1=\mathbf{k}_{t}(\alpha)$ and $\mathbf{a}_2=\mathbf{k}_{t}(\alpha-t)$ for some $\alpha$.
  \item[\rm(d)] $\mathbf{a}_1$ and $\mathbf{a}_2$ have a common eigenvector if and only if ${\rm tr}(\mathbf{a}_1\mathbf{a}_2)\in\{2,t^2-2\}$.
\end{enumerate}
\end{lem}

\begin{proof}
Suppose
$\mathbf{a}_i=\left(\begin{array}{cc} a_i & b_i \\ c_i & d_i \end{array}\right)$, with $a_id_i-b_ic_i=1$ and $a_i+d_i=t.$

(a) It follows from $\mathbf{a}_2=\mathbf{a}_1^{-1}\mathbf{d}(\lambda)$ that $a_2=\lambda d_1$ and $d_2=\lambda^{-1}a_1$, so
\begin{align*}
a_i=\frac{\lambda t}{\lambda+1}, \qquad d_i=\frac{t}{\lambda+1}, \qquad  b_ic_i=\frac{t^2-\lambda-\lambda^{-1}-2}{\lambda+\lambda^{-1}+2}\ne 0.
\end{align*}
Hence $\mathbf{a}_1=\mathbf{h}_{t}^{\lambda}(\mu)$ and $\mathbf{a}_2=\mathbf{h}_{t}^{\lambda}(-\lambda^{-1}\mu)$, with $\mu=(\lambda+1)b_1$.

(b) It follows from $\mathbf{a}_2=\mathbf{a}_1^{-1}\mathbf{p}$ that $a_2=d_1$ and $d_2=a_1-c_1$, so $c_1=0$. Then $c_2=0$, too.

(c) It follows from $\mathbf{a}_2=-\mathbf{a}_1^{-1}\mathbf{p}$ that $a_2=-d_1$ and $d_2=c_1-a_1$, so $c_1=2t$. Hence $\mathbf{a}_1=\mathbf{k}_{t}(\alpha)$ and $\mathbf{a}_2=\mathbf{k}_{t}(\lambda-t)$, with $\alpha=a_1-t/2$.

(d) ($\Rightarrow$) is obvious.

For ($\Leftarrow$), suppose ${\rm tr}(\mathbf{a}_1\mathbf{a}_2)\in\{2,t^2-2\}$. Replacing $\mathbf{a}_1$ by $\mathbf{a}_1^{-1}$ if necessary, we may assume ${\rm tr}(\mathbf{a}_1\mathbf{a}_2)=2$.
If $\mathbf{a}_1\mathbf{a}_2\ne\mathbf{e}$, then up to conjugacy we may assume $\mathbf{a}_1\mathbf{a}_2=\mathbf{p}$. By (b), $\mathbf{a}_1$ shares an eigenvector with $\mathbf{a}_2$.
\end{proof}

\begin{rmk}\label{rmk:generic}
\rm In the generic case, whenever ${\rm tr}(\mathbf{a}_1\mathbf{a}_2)=\lambda+\lambda^{-1}$, up to conjugacy we may assume $\mathbf{a}_1\mathbf{a}_2=\mathbf{d}(\lambda)$, so that $\mathbf{a}_1=\mathbf{h}_{t}^{\lambda}(\mu)$ and $\mathbf{a}_2=\mathbf{h}_{t}^{\lambda}(-\lambda^{-1}\mu)$ for some $\mu\ne 0$. This will simplify computations in many situations.
\end{rmk}

For $r=\eta+\eta^{-1}$ and $n\in\mathbb{Z}$, put
\begin{align*}
\omega_n(r)&=\begin{cases} (\eta^n-\eta^{-n})/(\eta-\eta^{-1}), &\eta\notin\{\pm 1\}, \\ n\eta^{n-1}, &\eta\in\{\pm 1\}, \end{cases} \\
\theta_n(r)&=\eta^n+\eta^{-n}.
\end{align*}
For any $\mathbf{z}\in{\rm SL}(2,\mathbb{C})$, by repeated applications of Hamilton-Cayley Theorem,
\begin{align}
\mathbf{z}^n=\omega_n({\rm tr}(\mathbf{z}))\mathbf{z}-\omega_{n-1}({\rm tr}(\mathbf{z}))\mathbf{e}.  \label{eq:power}
\end{align}

\section{Representations of rational tangles}

\begin{figure} [h]
  \centering
  \includegraphics[width=3cm]{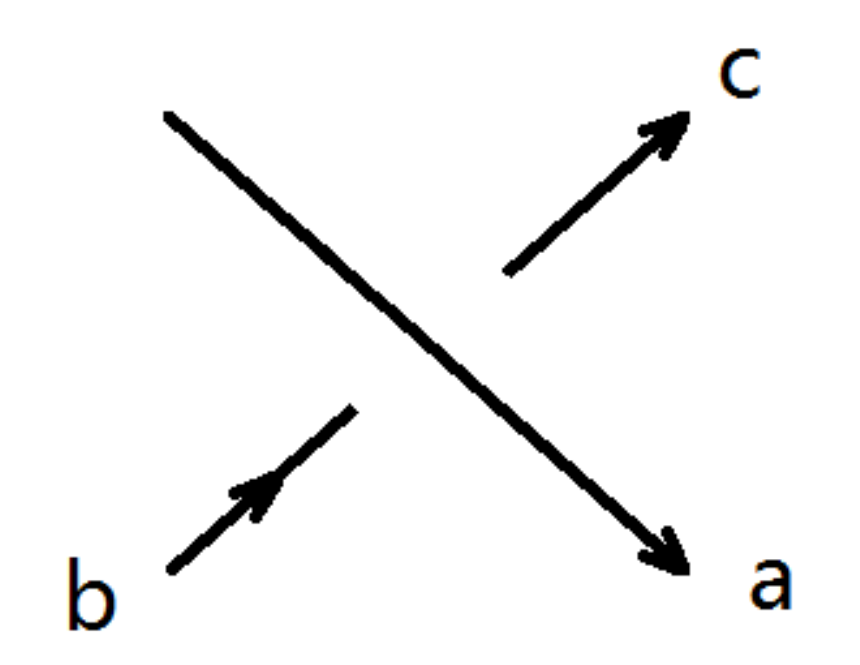}\\
  \caption{A representation satisfies $\rho(c)=\rho(a)\rho(b)\rho(a)^{-1}$ for each crossing.} \label{fig:crossing}
\end{figure}

Given a tangle $T$, let $\mathfrak{D}(T)$ denote the set of directed arcs of $T$; each arc gives rise to two directed arcs. If the directed arcs $a,b$ belong to the same component of $T$ (regarded as an embedded 1-manifold), then denote $[a]=[b]$ (resp. $[a]=-[b]$) if their directions are the same (resp. the opposite).
By a {\it representation} of $T$, we mean a map $\rho:\mathfrak{D}(T)\to {\rm SL}(2,\mathbb{C})$ such that $\rho(a^{-1})=\rho(a)^{-1}$ for each $a\in\mathfrak{D}(T)$, and $\rho(c)=\rho(a)\rho(b)\rho(a)^{-1}$ for each crossing as illustrated in Figure \ref{fig:crossing}.
To present such a representation, it is sufficient to give each arc a direction and label an element of ${\rm SL}(2,\mathbb{C})$ beside it.
We say that $\rho$ is {\it trace-$t$} if ${\rm tr}(\rho(a))=t$ for each $a\in\mathfrak{D}(T)$.

In virtue of Wirtinger presentation, a representation $\rho$ of $T$ can be identified with a true representation $\rho:\pi_1(B-T)\to{\rm SL}(2,\mathbb{C})$, where $B$ is a 3-ball containing $T$ such that $\partial B\cap T$ exactly consists of the end points of $T$.

Given a representation $\rho$ of $T\in\mathcal{T}_2^2$, let
\begin{align*}
\mathbf{x}^{{\rm nw}}=\rho(T^{{\rm nw}}), \qquad &\mathbf{x}^{{\rm ne}}=\rho(T^{{\rm ne}}), \qquad
\mathbf{x}^{{\rm sw}}=\rho(T^{{\rm sw}}), \qquad \mathbf{x}^{{\rm se}}=\rho(T^{{\rm se}}),
\end{align*}
where $T^{{\rm nw}}$, $T^{{\rm ne}}$, $T^{{\rm sw}}$, $T^{{\rm se}}$ are respectively the northwest, northeast, southwest, southeast ends of $T$, all directed outward.
Clearly,
$\mathbf{x}^{{\rm nw}}\mathbf{x}^{{\rm ne}}\mathbf{x}^{{\rm se}}\mathbf{x}^{{\rm sw}}=\mathbf{e}$,
so that $\mathbf{x}^{{\rm nw}}\mathbf{x}^{{\rm ne}}=(\mathbf{x}^{{\rm sw}})^{-1}(\mathbf{x}^{{\rm se}})^{-1}$.
Let
$$z={\rm tr}(\mathbf{x}^{{\rm nw}}\mathbf{x}^{{\rm ne}})-2, \qquad \dot{z}={\rm tr}(\mathbf{x}^{{\rm nw}}\mathbf{x}^{{\rm sw}})-2, \qquad
\grave{z}={\rm tr}(\mathbf{x}^{{\rm nw}}\mathbf{x}^{{\rm se}})-2.$$

\begin{figure}[h]
  \centering
  \includegraphics[width=10cm]{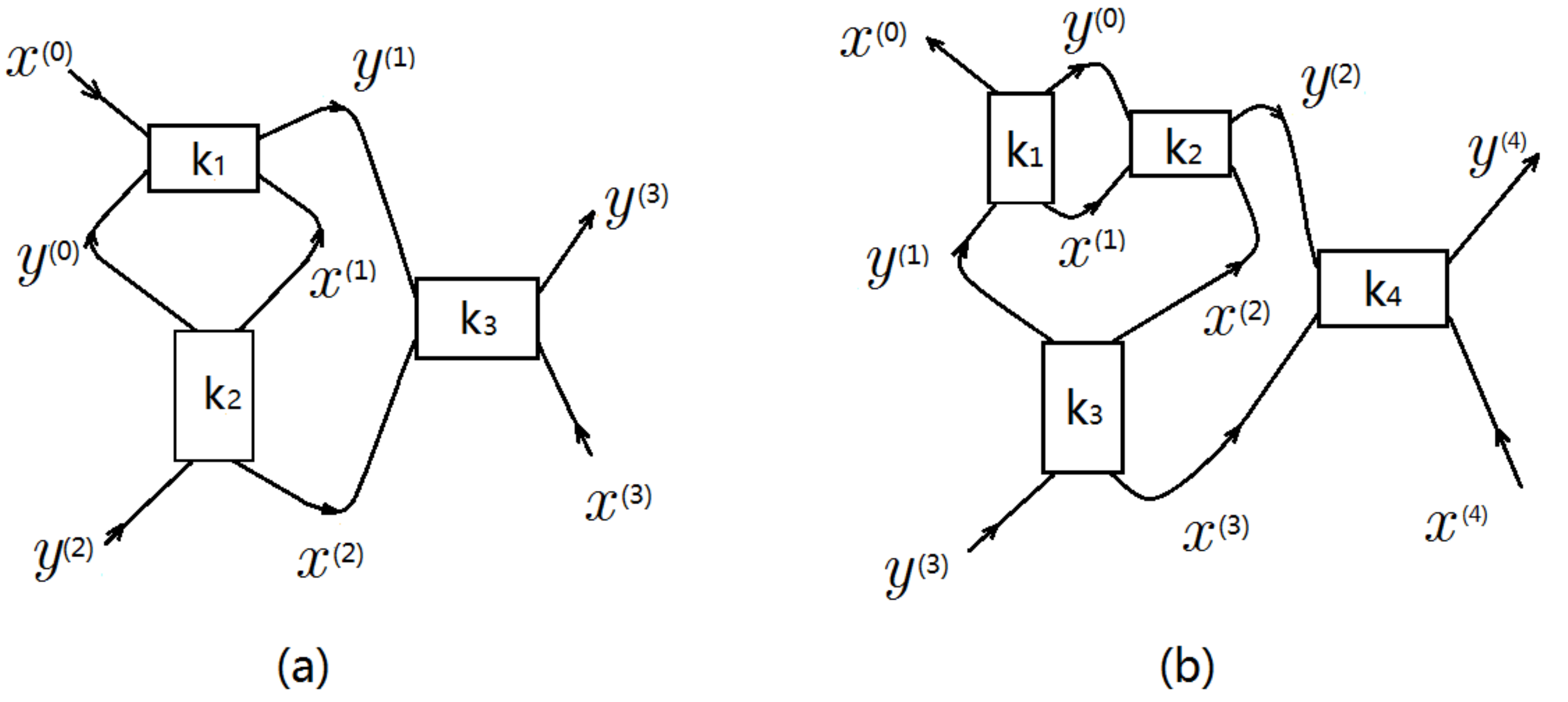}
  \caption{A representation of $[[k_1],\ldots,[k_s]]$: (a) $s$ is odd; (b) $s$ is even.} \label{fig:rational}
\end{figure}

\subsection{Some useful formulas}

From now on till the end of this section, let $T=[[k_1],\ldots,[k_s]]$ be a rational tangle, some of whose directed arcs are illustrated in Figure \ref{fig:rational}. 
We use a normal letter $a$ to denote a directed arc, and use boldface one $\mathbf{a}$ to denote the value $\rho$ takes at $a$. Let $x=T^{\rm nw}=(x^{(0)})^{(-1)^s}$ and $y=(y^{(0)})^{\pm 1}$, where the sign is chosen to make $[y]=[T^{\rm ne}]$ or $[y]=[T^{\rm sw}]$ if $\pm[x]=[T^{\rm ne}]$. Suppose $\rho$ is a representation of $T$.
Call $(\mathbf{x},\mathbf{y})$ the {\it generating pair} for $\rho$. Let $t={\rm tr}(\mathbf{x})={\rm tr}(\mathbf{y})$ and $r={\rm tr}(\mathbf{x}\mathbf{y})$.

\begin{lem} \label{lem:trivial}
$\mathbf{y}$ can be written as a word in $\mathbf{x}^{\rm ne}$, $\mathbf{x}^{\rm sw}$, $\mathbf{x}^{\rm se}$.
Consequently, $\rho$ is determined by $\mathbf{x}^{\rm nw}$, $\mathbf{x}^{\rm ne}$, $\mathbf{x}^{\rm sw}$, $\mathbf{x}^{\rm se}$.
\end{lem}

\begin{proof}
By induction on $i$, starting with $i=0$, it is easy to show that $\mathbf{x}^{(s-i)},\mathbf{y}^{(s-i)}$ can be written as words in $\mathbf{e},\mathbf{x}^{\rm ne},\mathbf{x}^{\rm sw},\mathbf{x}^{\rm se}$.
\end{proof}

Using (\ref{eq:power}) and the identity
$$\mathbf{y}\mathbf{x}=(r-t^2)\mathbf{e}+t(\mathbf{x}+\mathbf{y})-\mathbf{x}\mathbf{y},$$
in principle we can recursively compute $\mathbf{x}^{(j)}, \mathbf{y}^{(j)}$, $j=1,\ldots,s$, writing each one as a linear combination of $\mathbf{e},\mathbf{x}, \mathbf{y}, \mathbf{x}\mathbf{y}$, with coefficients in $\mathbb{Z}[t,r]$, and then obtain $z,\dot{z},\grave{z}\in\mathbb{Z}[t,r]$.
However, such computations are usually tedious, and they are unnecessary for our purpose. Here we take a better approach to derive the formulas for $z,\dot{z},\grave{z}$.

For $\gamma,\nu\ne 0$ and $n\in\mathbb{Z}$, introduce
\begin{align*}
\psi_\gamma^\nu(n)=\frac{2t^2+(2-\nu-\nu^{-1}-t^2)(\nu^n\gamma+\nu^{-n}\gamma^{-1})}{2-\nu-\nu^{-1}}.
\end{align*}

Let
$$-\nu_{j}-\nu_{j}^{-1}=e_j={\rm tr}\big(\mathbf{x}^{(-1)^j}\mathbf{y}^{(j)}\big), \qquad 0\le j\le s-1.$$

In the generic case, up to conjugacy we may assume $\mathbf{x}^{(0)}=\mathbf{h}_t^{-\nu_0}(1)$, $\mathbf{y}^{(0)}=\mathbf{h}_t^{-\nu_0}(\nu_0^{-1})$.
Then $\mathbf{y}^{(1)}=\mathbf{h}_t^{-\nu_0}(\nu_0^{k_1})$, $\mathbf{x}^{(1)}=\mathbf{h}_t^{-\nu_0}(\nu_0^{k_1-1})^{-1}$, and by (\ref{eq:tr-h}),
$$e_1=\psi_{1}^{\nu_0}(k_1), \qquad {\rm tr}\big(\mathbf{x}^{(0)}\mathbf{x}^{(1)}\big)=\psi_{1}^{\nu_0}(k_1-1).$$

Up to another conjugacy, we may assume $(\mathbf{x}^{(0)})^{-1}=\mathbf{h}_t^{-\nu_1}(1)$, $\mathbf{y}^{(1)}=\mathbf{h}_t^{-\nu_1}(\nu_1^{-1})$, and moreover, $\mathbf{y}^{(0)}=\mathbf{h}_t^{-\nu_1}(\gamma_1)$, $\mathbf{x}^{(1)}=\mathbf{h}_t^{-\nu_1}(\nu_1^{-1}\gamma_1)$ for some $\gamma_1$.
Then
\begin{align*}
e_0=\psi_{\gamma_1}^{\nu_1}(0), \qquad  {\rm tr}\big(\mathbf{x}^{(0)}\mathbf{x}^{(1)}\big)=\psi_{\gamma_1}^{\nu_1}(-1).
\end{align*}
From $\psi_{\gamma_1}^{\nu_1}(0)=e_0$ and $\psi_{\gamma_1}^{\nu_1}(-1)={\rm tr}\big(\mathbf{x}^{(0)}\mathbf{x}^{(1)}\big)=\psi_{1}^{\nu_0}(k_1-1)$ we obtain
$$\Big(1-\frac{t^2}{2+e_1}\Big)\gamma_1^{\pm 1}=\frac{\nu_1^{\pm 1}e_0-\psi_{1}^{\nu_0}(k_1-1)}{\nu_1^{\pm1}-\nu_1^{\mp1}}+\frac{2t^2}{(\nu_1^{\pm1}-\nu_1^{\mp1})(1-\nu_1^{\mp1})}.$$
Consequently, for each $n\in\mathbb{Z}$,
\begin{align*}
\psi_{\gamma_1}^{\nu_1}(n)&=\Big(1-\frac{t^2}{2+e_1}\Big)\gamma_1\nu_1^n+\Big(1-\frac{t^2}{2+e_1}\Big)\gamma_1^{-1}\nu_1^{-n}+\frac{2t^2}{2+e_1} \\
&=e_0\omega_{n+1}(-e_1)-\psi_{1}^{\nu_0}(k_1-1)\omega_n(-e_1)+\frac{2t^2}{2+e_1}(1-\omega_{n+1}(-e_1)+\omega_n(-e_1)).
\end{align*}
Moreover, $\mathbf{y}^{(2)}=\mathbf{h}_t^{-\nu_1}(\nu_1^{k_2}\gamma_1)$, $\mathbf{x}^{(2)}=\mathbf{h}_t^{-\nu_1}(\nu_1^{k_2-1}\gamma_1)^{-1}$, and
$$e_2=\psi_{\gamma_1}^{\nu_1}(k_2), \qquad {\rm tr}\big((\mathbf{x}^{(0)})^{-1}\mathbf{x}^{(2)}\big)=\psi_{\gamma_1}^{\nu_1}(k_2-1).$$

Going on, up to yet another conjugacy, we may assume
$$\mathbf{x}^{(0)}=\mathbf{h}_t^{-\nu_2}(1),  \ \ \   \mathbf{y}^{(2)}=\mathbf{h}_t^{-\nu_2}(\nu_2^{-1}), \ \ \ \mathbf{y}^{(1)}=\mathbf{h}_t^{-\nu_2}(\gamma_2), \ \ \ \mathbf{x}^{(2)}=\mathbf{h}_t^{-\nu_2}(\nu_2^{-1}\gamma_2),$$
for some $\gamma_2$. Then
\begin{align*}
e_1=\psi_{\gamma_2}^{\nu_2}(0), \qquad
{\rm tr}\big((\mathbf{x}^{(0)})^{-1}\mathbf{x}^{(2)}\big)=\psi_{\gamma_2}^{\nu_2}(-1),
\end{align*}
which, similarly as above, imply that for each $n\in\mathbb{Z}$,
\begin{align*}
\psi_{\gamma_2}^{\nu_2}(n)=e_1\omega_{n+1}(-e_2)-\psi_{\gamma_1}^{\nu_1}(k_2-1)\omega_n(-e_2)+\frac{2t^2}{2+e_2}(1-\omega_{n+1}(-e_2)+\omega_n(-e_2)).
\end{align*}

In general, there exist $\gamma_0=1$ and $\gamma_1,\ldots,\gamma_{s-1}\in\mathbb{C}^\ast$ such that, with $\psi_{\gamma_j}^{\nu_j}(n)$ abbreviated to $\psi_j(n)$,
\begin{align}
\psi_0(n)&=\theta_n(-e_0)+\frac{2-\theta_n(-e_0)}{2+e_0}t^2,    \label{eq:useful-0}  \\
e_{j}&=\psi_{j-1}(k_j), \qquad 1\le j\le s-1,  \label{eq:useful-1}  \\
\psi_j(n)&=e_{j-1}\omega_{n+1}(-e_j)-\psi_{j-1}(k_j-1)\omega_n(-e_j)+\frac{2t^2(1-\omega_{n+1}(-e_j)+\omega_n(-e_j))}{2+e_j}.  \label{eq:useful-2}
\end{align}
So the $\psi_j(n)$'s can be recursively computed as a polynomial in $t$ and $r$, without explicitly referring to $\gamma_j$.
In particular, we can obtain the expressions for
$$\dot{z}=e_{s-1}-2, \qquad \grave{z}=\psi_{s-1}(k_s-1)-2, \qquad z=\psi_{s-1}(k_s)-2.$$
Note that these equations are polynomial equations. Since they are available under some Zariski open conditions, actually they always hold.

\begin{rmk}\label{rmk:trace-free}
\rm When $t=0$, (\ref{eq:useful-0})--(\ref{eq:useful-2}) can be much simplified and can lead to closed formulas for $z,\dot{z},\grave{z}$ in terms of $\nu_0$.
\end{rmk}

Write $t=\kappa+\kappa^{-1}$.
If $T$ is odd, then $N(T)$ is a knot. We assert that
$$z=(r-2)\phi_{N(T)}(\kappa,2-r)^2,$$
where 
$\phi_{N(T)}(\kappa,u)$ is the {\it Riley polynomial} of $N(T)$, as introduced in \cite{Ri84}.

It can be shown that $x^{\rm ne}=vyv^{-1}$ for some word
$$v=(y^{-1})^{\epsilon_1}x^{\epsilon_2}\cdots(y^{-1})^{\epsilon_{2h-1}}x^{\epsilon_{2h}}, \qquad \text{with} \quad \epsilon_i=\epsilon_{2h+1-i}\in\{\pm 1\}.$$
For a nonabelian representation of $N(T)$, up to conjugacy one may assume $\mathbf{x}=\mathbf{u}^-_\kappa(u)$, $\mathbf{y}^{-1}=\mathbf{u}^+_\kappa(1)$, so that $u=2-r$. Suppose $\mathbf{v}=(v_{ij})_{2\times 2}$. Taking a square root $\tilde{u}$ of $u$ (and using prime to denote the transpose), we deduce from
$$\mathbf{v}'=(\mathbf{x}')^{\epsilon_1}\cdots(\mathbf{y}')^{-\epsilon_{2h}}=\mathbf{d}(\tilde{u})\mathbf{v}\mathbf{d}(\tilde{u})^{-1}$$
that $v_{21}=uv_{12}$. By definition (slightly different from the original one given on Page 7 of \cite{Ri84}), $\phi_{N(T)}(\kappa,u)=v_{11}+(\kappa^{-1}-\kappa)v_{12}$.
Then
$$z={\rm tr}(\mathbf{x}\mathbf{v}\mathbf{y}\mathbf{v}^{-1})-2=-u(v_{11}+(\kappa^{-1}-\kappa)v_{12})^2,$$
proving the assertion.

If $T$ is even, then $N(T)$ has two components, and
\begin{align*}
z=(r-2)(r+2-t^2)\phi^\iota_{N(T)}(\kappa,r-\kappa^{1+\iota}-\kappa^{-1-\iota})^2, 
\end{align*}
where $\iota\in\{\pm 1\}$ and $\phi^\iota_{N(T)}(\kappa,u)$ is defined as follows.

We can show that $x^{\rm ne}=vx^{-1}v^{-1}$ for some word $v=y^{\epsilon_1}x^{\epsilon_2}\cdots x^{\epsilon_{2h-2}}y^{\epsilon_{2h-1}}$ with $\epsilon_i=\epsilon_{2h-i}\in\{\pm1\}$ for all $i$.
Suppose
$$\mathbf{x}=\mathbf{u}^-(\kappa,u), \qquad \mathbf{y}=\mathbf{u}^+(\kappa^\iota,1), \qquad  u=r-\kappa^{1+\iota}-\kappa^{-1-\iota}.$$
If $\mathbf{v}=(v_{ij})_{2\times 2}$, then put $\phi^\iota_{N(T)}(\kappa,u)=v_{12}$.

Take
$$\mathbf{a}=\left(\begin{array}{cc} (\kappa^\iota-\kappa^{-\iota})u & u \\ u & \kappa^{-1}-\kappa \end{array}\right).$$
Then $\mathbf{a}\mathbf{x}=\mathbf{x}'\mathbf{a}$ and $\mathbf{a}\mathbf{y}=\mathbf{y}'\mathbf{a}$. Hence $\mathbf{a}\mathbf{v}=\mathbf{v}'\mathbf{a}$, which implies
$$(\kappa-\kappa^{-1})v_{21}=(v_{11}-v_{22}+(\kappa^{-\iota}-\kappa^\iota)v_{12})u.$$
Due to this, $\mathbf{x}\mathbf{v}=\mathbf{v}\mathbf{x}$ if and only if $v_{12}=0$. One has
$$z={\rm tr}(\mathbf{x}\mathbf{v}\mathbf{x}^{-1}\mathbf{v}^{-1})-2=(r-\kappa^{1+\iota}-\kappa^{-1-\iota})(r-\kappa^{1-\iota}-\kappa^{\iota-1})v_{12}^2.$$

\begin{rmk} \label{rmk:reducible1}
\rm When $T$ is odd, non-abelian representations are given by $\phi_{N(T)}=0$. In particular, $\phi_{N(T)}(\kappa,0)=\pm\kappa^{n}\Delta_{N(T)}(\kappa^2)$ for some $n\in\mathbb{Z}$, where $\Delta_{N(T)}$ is the Alexander polynomial of $N(T)$.

When $T$ is even, non-abelian representations are given by $\phi^\iota_{N(T)}=0$. In particular, $\phi^{\pm}_{N(T)}(\kappa,0)=\kappa^n\Delta_{N(T)}(\kappa^2,\kappa^{\pm2})$ for some $n\in\mathbb{Z}$, where $\Delta_{N(T)}$ is the two-variable Alexander polynomial of $N(T)$.
\end{rmk}

\begin{exmp}\label{exmp:odd-twist}
\rm Let $T=[2h+1]$. Clearly, $e_0=t^2-r$. Then
\begin{align*}
\dot{z}&=t^2-r-2, \\
\grave{z}&=(r-2)(r+2-t^2)\omega_h(t^2-r)^2, \\
z&=(r-2)\big(\omega_{h+1}(t^2-r)-\omega_h(t^2-r)\big)^2.
\end{align*}

Putting $t=\kappa+\kappa^{-1}$ and using the identity
\begin{align*}
\omega_{k+1}(a)^2+\omega_{k}(a)^2=1+a\omega_{k+1}(a)\omega_{k}(a),
\end{align*}
we can rewrite
\begin{align*}
z=t^2-4+(2-e_0)\big(\omega_{h+1}(e_0)-\kappa^2\omega_h(e_0)\big)\big(\omega_{h+1}(e_0)-\kappa^{-2}\omega_h(e_0)\big). 
\end{align*}
\end{exmp}

\begin{exmp}\label{exmp:double-twist}
\rm Let $T=[[k_1],[k_2]]$.
\begin{align*}
\dot{z}&=e_1-2=\psi_{0}(k_1)-2=\frac{(2+e_0-t^2)(\theta_{k_1}(-e_0)-2)}{2+e_0}, \\
\grave{z}&=\psi_{1}(k_2-1)-2, \\
z&=\psi_{1}(k_2)-2,
\end{align*}
with
\begin{align*}
\psi_1(n)=\ &r\omega_{n+1}(-e_1)-\Big(\frac{2-\theta_{k_1-1}(-e_0)}{2+e_0}t^2+\theta_{k_1-1}(-e_0)\Big)\omega_n(-e_1) \\
&+\frac{2(1-\omega_{n+1}(-e_1)+\omega_n(-e_1))}{2+e_1}t^2.
\end{align*}
In particular,
\begin{itemize}
  \item when $k_1=k_2=2$, we have $e_0=r$, $e_1=r^2-2+(2-r)t^2$,
        \begin{align*}
        \grave{z}&=-re_1-(t^2-r)+2t^2-2=(t^2-r-2)(r-1)^2, \\
        z&=r(e_1^2-e_1-1)+t^2(2-e_1)-2=(r-2)(r^2+(r-1)(1-t^2))^2.
        \end{align*}
  \item when $k_1=3,k_2=2$, we have $e_0=r$, $e_1=3r-r^3+t^2(r-1)^2$,
        \begin{align*}
        \grave{z}&=-re_1-((2-r)t^2+r^2-2)+2t^2-2=r^2(r-2)(r+2-t^2), \\
        z&=r(e_1^2-1)+(r^2-t^2r-2)e_1+2t^2-2 \\
        &=(r-2)(r^3+(1-t^2)r^2+(t^2-2)r-1)^2.
        \end{align*}
\end{itemize}
\end{exmp}

\subsection{A class of reducible representations}

When $\rho$ is nonabelian reducible with $z=t^2-4\ne 0,-4$ (so that $T$ is odd), we are no longer able to know $\rho$ well enough merely from traces, and need to look into it more closely.

Up to conjugacy, we can assume $\mathbf{x}=\mathbf{u}^+_{\kappa}(0)=\mathbf{d}(\kappa)$ and $\mathbf{y}=\mathbf{u}^+_{\kappa}(1)$, where $\kappa+\kappa^{-1}=t$.
Let $\mathbf{a}\lrcorner\mathbf{b}=\mathbf{a}\mathbf{b}\mathbf{a}^{-1}$. Using
\begin{align*}
\mathbf{u}^+_\mu(a)^k&=\mathbf{u}^+_{\mu^k}(\gamma_k(\mu+\mu^{-1})a), \\
\mathbf{u}^+_\mu(a)\lrcorner\mathbf{u}^+_\nu(b)&=\mathbf{u}^+_\nu\big((\nu^{-1}-\nu)\mu a+\mu^2b\big),
\end{align*}
we can recursively compute $\mathbf{y}^{(j)}$, so as to obtain $\mathbf{x}^{\rm ne}=\mathbf{y}^{(s)}$, $\mathbf{x}^{\rm sw}=(\mathbf{y}^{(s-1)})^{-1}$.
Such computations are much easier than those for irreducible representations.
Let $\vartheta_T^{\rm ne}(\kappa)$ and $\vartheta_T^{\rm sw}(\kappa)$ respectively denote the upper-right entry of $\mathbf{y}^{(s)}$ and $\mathbf{y}^{(s-1)}$.

\begin{exmp}\label{exmp:red-1}
\rm Let $T=[2h+1]$. 
Trivially, $\vartheta_T^{\rm sw}(\kappa)=1$. Since
\begin{align*}
\mathbf{y}^{(1)}&=(\mathbf{x}^{-1}\mathbf{y})^{h+1}\lrcorner\mathbf{y}=\mathbf{u}^+_{1}(\kappa^{-1})^{h+1}\lrcorner\mathbf{u}^+_{\kappa}(1)
=\mathbf{u}^+_{\kappa}((h+1)\kappa^{-2}-h),
\end{align*}
we have $\vartheta_T^{\rm ne}(\kappa)=(h+1)\kappa^{-2}-h$.
\end{exmp}

\begin{exmp}\label{exmp:red-2}
\rm Let $T=[[2h_1],[2h_2]]$. 
Clearly,
\begin{align*}
\mathbf{y}^{(1)}&=(\mathbf{x}\mathbf{y})^{h_1}\lrcorner\mathbf{x}=\mathbf{u}^+_{\kappa^2}(\kappa)^{h_1}\lrcorner\mathbf{u}^+_{\kappa}(0)
=\mathbf{u}^+_{\kappa}\big((1-\kappa^{4h_1})/(1+\kappa^{-2})\big), \\
\mathbf{y}^{(2)}&=(\mathbf{x}^{-1}\mathbf{y}^{(1)})^{h_2}\lrcorner\mathbf{y}
=\mathbf{u}^+_{1}\big(h_2(1-\kappa^{4h_1})/(\kappa+\kappa^{-1})\big)\lrcorner\mathbf{u}^+_{\kappa}(1) \\
&=\mathbf{u}^+_{\kappa}\Big(1+h_2\frac{(\kappa^2-1)(\kappa^{4h_1}-1)}{\kappa^{2}+1}\Big).
\end{align*}
Hence
$$\vartheta_T^{\rm ne}(\kappa)=1+h_2\frac{(\kappa^2-1)(\kappa^{4h_1}-1)}{\kappa^{2}+1}, \qquad
\vartheta_T^{\rm sw}(\kappa)=\frac{1-\kappa^{4h_1}}{1+\kappa^{-2}}.$$
\end{exmp}

\begin{exmp}\label{exmp:red-3}
\rm For $T=[[2h_1+1],[2h_2]]$, 
\begin{align*}
\mathbf{y}^{(1)}&=(\mathbf{x}\mathbf{y})^{h_1+1}\lrcorner\mathbf{y}=\mathbf{u}^+_{\kappa^2}(\kappa)^{h_1+1}\lrcorner\mathbf{u}^+_{\kappa}(1)
=\mathbf{u}^+_{\kappa}\big((1+\kappa^{4h_1+2})/(1+\kappa^{-2})\big), \\
\mathbf{y}^{(2)}&=(\mathbf{x}^{-1}\mathbf{y}^{(1)})^{h_2}\lrcorner\mathbf{y}
=\mathbf{u}^+_{1}\big(h_2(1+\kappa^{4h_1+2})/(\kappa+\kappa^{-1})\big)\lrcorner\mathbf{u}^+_{\kappa}(1) \\
&=\mathbf{u}^+_{\kappa}\Big(1+h_2\frac{(1-\kappa^2)(1+\kappa^{4h_1+2})}{1+\kappa^{2}}\Big).
\end{align*}
Hence
$$\vartheta_T^{\rm ne}(\kappa)=1+h_2\frac{(1-\kappa^2)(1+\kappa^{4h_1+2})}{1+\kappa^{2}}, \qquad
\vartheta_T^{\rm sw}(\kappa)=\frac{1+\kappa^{4h_1+2}}{1+\kappa^{-2}}.$$
\end{exmp}

The following will be applied in Section \ref{sec:X'}.
\begin{lem}\label{lem:reducible2}
Suppose $\mathbf{x}=\mathbf{u}^\varepsilon_\kappa(\alpha)$, $\mathbf{y}=\mathbf{u}^\varepsilon_{\kappa}(\beta)$, with $\varepsilon\in\{+,-\}$, $\alpha\ne\beta$. Then (with $\pm$ identified with $\pm1$)
\begin{enumerate}
  \item[\rm(i)] $\mathbf{x}^{\rm ne}=\mathbf{x}$ if and only if $\vartheta_T^{\rm ne}(\kappa^\varepsilon)=0$;
  \item[\rm(ii)] $\mathbf{x}^{\rm se}=\mathbf{x}^{\rm sw}$ if and only if $\vartheta^{\rm ne}_T(\kappa^\varepsilon)=(1+\kappa^{-2\varepsilon})\vartheta^{\rm sw}_T(\kappa^\varepsilon)$.
\end{enumerate}
\end{lem}

\begin{proof}
Using the identity $\mathbf{w}\lrcorner\mathbf{u}^-_\kappa(\xi)=\mathbf{u}^+_{1/\kappa}(-\xi)$ if necessary, we may assume $\varepsilon=+$.

Recall that in the special case $\mathbf{x}=\mathbf{u}^+_\kappa(0)$, $\mathbf{y}=\mathbf{u}^+_{\kappa}(1)$, we have
$$\mathbf{x}^{\rm ne}=\mathbf{u}^+_{\kappa}(\vartheta_T^{\rm ne}(\kappa)), \qquad  \mathbf{x}^{\rm sw}=\mathbf{u}^+_{\kappa}(\vartheta_T^{\rm sw}(\kappa))^{-1}.$$
Take $\eta,\gamma$ with $\eta^2=\beta-\alpha$ and $\eta(\kappa^{-1}-\kappa)\gamma=\alpha$. We have
\begin{align*}
\mathbf{u}^+_\eta(\gamma)\mathbf{u}^+_\kappa(0)\mathbf{u}^+_\eta(\gamma)^{-1}=\mathbf{u}^+_\kappa(\alpha),  \qquad
\mathbf{u}^+_\eta(\gamma)\mathbf{u}^+_\kappa(1)\mathbf{u}^+_\eta(\gamma)^{-1}=\mathbf{u}^+_\kappa(\beta).
\end{align*}
Thus,
\begin{align*}
\mathbf{x}^{\rm ne}&=\mathbf{u}^+_{\kappa}\big((1-\vartheta_T^{\rm ne}(\kappa))\alpha+\vartheta_T^{\rm ne}(\kappa)\beta\big),    \\
\mathbf{x}^{\rm sw}&=\mathbf{u}^+_{\kappa}\big((1-\vartheta_T^{\rm sw}(\kappa))\alpha+\vartheta_T^{\rm sw}(\kappa)\beta\big)^{-1}.
\end{align*}
Then (i) is clear, and (ii) is seen from $\mathbf{x}^{\rm se}=(\mathbf{x}^{\rm sw}\mathbf{x}\mathbf{x}^{\rm ne})^{-1}$.
\end{proof}

\section{Irreducible trace-$t$ representations of Montesinos knots with $t\ne 0$}

The case when $t=0$ has been well-understood \cite{Ch18-1}. We do not recall it here, but point out that the result of \cite{Ch18-1} can be quickly recovered, in view of Remark \ref{rmk:trace-free}.

Throughout this section, assume $t\ne 0$.

Consider a Montesinos knot $K=M(p_1/q_1,\ldots,p_m/q_m)$. Let $T_i=[p_i/q_i]$. Adopt the cyclical notation for subscript indices, so that $T_{m+1}$ means $T_1$, etc.

Let $(\mathbf{x}_i,\mathbf{y}_i)$ be the generating pair for $\rho|_{T_i}$ (whose meaning is self-evident), and let $r_i={\rm tr}(\mathbf{x}_i\mathbf{y}_i)$.
Let $\mathbf{x}_i^{\rm ne}=\rho(T_i^{\rm ne})$, etc. Clearly, $\mathbf{x}_i^{\rm sw}=(\mathbf{x}_{i+1}^{\rm nw})^{-1}$ and
$\mathbf{x}_i^{\rm se}=(\mathbf{x}_{i+1}^{\rm ne})^{-1}$.
Let
$$z_i-2={\rm tr}(\mathbf{x}_i^{\rm nw}\mathbf{x}_i^{\rm ne}), \qquad
\dot{z}_i-2={\rm tr}(\mathbf{x}_i^{\rm nw}\mathbf{x}_i^{\rm sw}), \qquad
\grave{z}_i-2={\rm tr}(\mathbf{x}_i^{\rm nw}\mathbf{x}_i^{\rm se}).$$

Let $\mathbf{g}=\mathbf{x}_1^{\rm nw}\mathbf{x}_1^{\rm ne}=\cdots=\mathbf{x}_m^{\rm nw}\mathbf{x}_m^{\rm ne}$.
Suppose ${\rm tr}(\mathbf{g})=\tau=\lambda+\lambda^{-1}$, so that $z_i=\tau-2$ for all $i$.

\begin{lem}\label{lem:elementary}
{\rm(a)} If $\tau=2$, then $\mathbf{g}=\mathbf{e}$.

{\rm(b)} $\mathbf{g}\ne-\mathbf{e}$.
\end{lem}
\begin{proof}
(a) If $\tau=2$ but $\mathbf{g}\ne\mathbf{e}$, then up to conjugacy we may assume $\mathbf{g}=\mathbf{p}$. By Lemma \ref{lem:key} (b), $\mathbf{x}_i^{\rm nw},\mathbf{x}_i^{\rm ne}$, $i=1,\ldots,m$, are all upper-triangular, so by Lemma \ref{lem:trivial}, $\rho$ would be reducible, contradicting the assumption.

(b) If $\mathbf{g}=-\mathbf{e}$, then $t={\rm tr}((\mathbf{x}_1^{\rm ne})^{-1})={\rm tr}(-\mathbf{x}_1^{\rm nw})=-t$, which is impossible under the hypothesis $t\ne 0$.
\end{proof}

Separately, there are three possibilities: (i) $\mathbf{g}=\mathbf{e}$, (ii) $\tau\notin\{2, t^2-2\}$, and (iii) $\tau=t^2-2\ne\pm 2$, which respectively give rise to $\mathcal{X}_1(K)$, $\mathcal{X}_2(K)$ and $\mathcal{X}'(K)$.

\subsection{$\mathbf{g}=\mathbf{e}$}

In this case, $\rho|_{T_i}$ induces a representation of $N(T_i)$, which is denoted by $\check{\rho}_i$. Clearly, for each $i\in\{1,\ldots,m\}$,
$$\check{\rho}_{i+1}(T_{i+1}^{\rm nw})=\check{\rho}_i(T_i^{\rm sw})^{-1}, \qquad \check{\rho}_{i+1}(T_{i+1}^{\rm ne})=\check{\rho}_i(T_i^{\rm se})^{-1}.$$
Conversely, given representations $\rho_i$ of $N(T_i)$ such that $\rho_{i+1}(T_{i+1}^{\rm nw})=\rho_i(T_i^{\rm sw})^{-1}$, $\rho_{i+1}(T_{i+1}^{\rm ne})=\rho_i(T_i^{\rm se})^{-1}$ for all $i$, there exists a unique representation $\rho$ of $K$ with $\check{\rho}_i=\rho_i$.
In this sense, we say that $\rho$ is the {\it union} of the $\rho_i$'s.


If $\rho|_{T_i}$ is abelian, then $\mathbf{x}_i^{\rm sw}=\mathbf{x}_i$ or $\mathbf{x}_i^{\rm sw}=\mathbf{x}_i^{-1}$. More, precisely, if $2\nmid p_i$, then $\mathbf{x}_i^{\rm sw}=\mathbf{x}_i^{(-1)^{q_i+1}}$;
if $2\mid p_i$, then both $\mathbf{x}_i^{\rm sw}=\mathbf{x}_i$ and $\mathbf{x}_i^{\rm sw}=\mathbf{x}_i^{-1}$ are possible.

If $\rho|_{T_i}$ is nonabelian reducible and $2\nmid p_i$, then (with $\kappa+\kappa^{-1}=t$ as always) $\Delta_{N(T_i)}(\kappa^2)=0$; if $\rho|_{T_i}$ is nonabelian reducible and $2\mid p_i$, then $\Delta_{N(T_i)}(\kappa^2,\kappa^2)=0$ or  $\Delta_{N(T_i)}(\kappa^2,\kappa^{-2})=0$.

Remember that $\rho$ is determined by $\mathbf{x}_1,\ldots,\mathbf{x}_m$.
Let
$$R_t=\{(r_1,\ldots,r_m)\colon z_i(t,r_i)=0,\ 1\le i\le m\}.$$
For each $\vec{r}=(r_1,\ldots,r_m)\in R_t$, let $\mathcal{T}(t,\vec{r})$ denote the space of $(\mathbf{x}_1,\ldots,\mathbf{x}_m)\in({\rm SL}(2,\mathbb{C})-\{\pm\mathbf{e}\})^m$ satisfying
\begin{align*}
{\rm tr}(\mathbf{x}_i)=t, \quad {\rm tr}(\mathbf{x}_i^{-1}\mathbf{x}_{i+1})=\dot{z}_i(t,r_i)+2, \qquad 1\le i\le m.
\end{align*}
Then trace-$t$ representations $\rho$ bijectively correspond to elements of $\bigcup_{\vec{r}\in R_t}\mathcal{T}(t,\vec{r}).$


Therefore, $\mathcal{X}_1(K)$ have high-dimensional components when $m\ge 4$.

Let $\overline{\mathbf{x}}_i=\mathbf{x}-(t/2)\mathbf{e}$. According to \cite{ABL18} Theorem 3.1, the character of the induced representation $F_m\to\pi(K)\stackrel{\rho}\to{\rm SL}(2,\mathbb{C})$ is determined by
\begin{align*}
\overline{t}_{i_1,i_2}:&={\rm tr}(\overline{\mathbf{x}}_{i_1}\overline{\mathbf{x}}_{i_2}), &&1\le i_1<i_2\le m,  \\
\overline{t}_{i_1,i_2,i_3}:&={\rm tr}(\overline{\mathbf{x}}_{i_1}\overline{\mathbf{x}}_{i_2}\overline{\mathbf{x}}_{i_3}),  &&1\le i_1<i_2<i_3\le m,
\end{align*}
subject to the following relations:
\begin{align}
2\overline{t}_{i_1,i_2,i_3}\overline{t}_{j_1,j_2,j_3}+\det\big((\overline{t}_{i_a,j_b})_{3\times 3}\big)=0, \label{eq:typeI}
\end{align}
for all $1\le i_1<i_2<i_3\le m$ and $1\le j_1<j_2<j_3\le m$;
\begin{align}
\overline{t}_{i,j_0}\overline{t}_{j_1,j_2,j_3}-\overline{t}_{i,j_1}\overline{t}_{j_0,j_2,j_3}
+\overline{t}_{i,j_2}\overline{t}_{j_0,j_1,j_3}-\overline{t}_{i,j_3}\overline{t}_{j_0,j_1,j_2}=0,  \label{eq:typeII}
\end{align}
for all $1\le i\le m$ and $1\le j_0<j_1<j_2<j_3\le m$.

We can describe $\mathcal{X}_1(K)$ using $t$, $\overline{t}_{i_1,i_2}$, $\overline{t}_{i_1,i_2,i_3}$, $r_i$, subject to
$$z_i=0, \qquad \frac{1}{2}t^2-\overline{t}_{i,i+1}=\dot{z}_i(t,r_i)+2,$$
and (\ref{eq:typeI}), (\ref{eq:typeII}), and conditions for irreducibility:
\begin{lem}
$\rho$ is irreducible if and only if one of the following holds: (i) $\rho|_{T_i}$ is irreducible for some $i$, which is equivalent to $r_{i}\notin\{2,t^2-2\}$; (ii) $r_{i}\in\{2,t^2-2\}$ for all $i$ and $\overline{t}_{i,i+1,i'}\ne 0$ for some $i,i'$.
\end{lem}
\begin{proof}
Notice that $(\mathbf{a},\mathbf{b},\mathbf{c})\mapsto{\rm tr}(\mathbf{a}\mathbf{b}\mathbf{c})$ defines a nontrivial alternating trilinear form on the 3-dimensional space $\{\mathbf{a}\in M\colon {\rm tr}(\mathbf{a})=0\}$.
Hence $\overline{\mathbf{x}}_{i_1},\overline{\mathbf{x}}_{i_2},\overline{\mathbf{x}}_{i_3}$ are linearly independent if and only if $\overline{t}_{i_1,i_2,i_3}\ne 0$.

If $\rho$ is irreducible and $\rho|_{T_{i}}$ is reducible for all $i$, then there exists $i$ with $\rho|_{T_{i}}$ nonabelian, so that $\overline{\mathbf{x}}_{i}\ne\pm\overline{\mathbf{x}}_{i+1}$. The irreducibility of $\rho$ implies the existence of $\overline{\mathbf{x}}_{i'}$ which is not a linear combination of $\overline{\mathbf{x}}_{i}$ and $\overline{\mathbf{x}}_{i+1}$. Then $\overline{t}_{i,i+1,i'}\ne 0$.
\end{proof}

\begin{rmk}\label{rmk:ensure-irr}
\rm If $K$ is odd, then $\rho$ has been guaranteed to be irreducible. Otherwise, up to conjugacy we may assume ${\rm Im}(\rho)\subset U$ and $\mathbf{x}_1=\mathbf{u}^+_\kappa(1)=(\mathbf{x}_1^{\rm ne})^{-1}$; on the other hand, given an orientation to $K$, due to $K\in\mathcal{M}_1$, the induced direction on $T_1^{\rm ne}$ is the same as that on $T_1^{\rm nw}$, so $\mathbf{x}_1^{\rm ne}$ can be conjugated to $\mathbf{x}_1$ via elements in ${\rm Im}(\rho)$. This is impossible.
\end{rmk}


\subsection{$\tau\ne 2,t^2-2$}\label{sec:X2}

When $\tau\ne\pm 2, t^2-2$, up to conjugacy we may assume $\mathbf{g}=\mathbf{d}(\lambda)$ with
\begin{align}
|\lambda|>1 \qquad \text{or} \quad |\lambda|=1, \ \arg\lambda\in(0,\pi]. \label{eq:lambda}
\end{align}
Then for each $i$ there exists $\mu_i\ne 0$ such that $\mathbf{x}_i^{\rm nw}=\mathbf{h}_t^{\lambda}(\mu_i)$, $\mathbf{x}_i^{\rm ne}=\mathbf{h}_t^{\lambda}(-\lambda^{-1}\mu_i)$. Applications of (\ref{eq:tr-h}) to $(\mu,\nu)=(\mu_i,\mu_{i+1})$ and $(\mu,\nu)=(\mu_i,-\lambda^{-1}\mu_{i+1})$ enable us to solve
$$\frac{\mu_i}{\mu_{i+1}}=\frac{(1+\lambda^{-1})\dot{z}_i+(1+\lambda)\grave{z}_i+2(\tau+2-t^2)}{(1-\lambda)(\tau+2-t^2)}.$$
Since $\prod_{i=1}^m(\mu_i/\mu_{i+1})=1$, we have
\begin{align}
\prod\limits_{i=1}^m\big((1+\lambda^{-1})\dot{z}_i+(1+\lambda)\grave{z}_i+2(\tau+2-t^2)\big)=(1-\lambda)^m(\tau+2-t^2)^m.  \label{eq:generic-0}
\end{align}

When $\tau=-2$, by Lemma \ref{lem:elementary} (b), up to conjugacy we may assume $\mathbf{g}=-\mathbf{p}$.
Then
$\mathbf{x}_i^{\rm nw}=\mathbf{k}_t(\alpha_i)$, $\mathbf{x}_i^{\rm ne}=\mathbf{k}_t(\alpha_i-t)$. Applications of (\ref{eq:tr-k}) to $(\alpha,\beta)=(\alpha_i,\alpha_{i+1})$ and $(\alpha,\beta)=(\alpha_i,\alpha_{i+1}-t)$ enable us to solve
$$2t(\alpha_{i+1}-\alpha_{i})=\dot{z}_i-\grave{z}_i+t^2.$$
Since $\sum_{i=1}^m(\alpha_{i+1}-\alpha_i)=0$, we have
\begin{align}
\sum\limits_{i=1}^m(\grave{z}_i-\dot{z}_i)=mt^2.   \label{eq:tau=-2}
\end{align}


We can unify (\ref{eq:generic-0}) and (\ref{eq:tau=-2}) as $H(\lambda,t;r_1,\ldots,r_m)=0$, with $H$ given by
\begin{align}
\frac{1}{\lambda+1}\left(\prod\limits_{i=1}^m\big((1+\lambda^{-1})\dot{z}_i+(1+\lambda)\grave{z}_i+2(\tau+2-t^2)\big)-(1-\lambda)^m(\tau+2-t^2)^m\right).  \label{eq:generic-1}
\end{align}
To see this, just note that $H$ is a polynomial in $\lambda+1$ with constant term
$$(-1)^m\Big(mt^2+\sum_{i=1}^m(\dot{z}_i-\grave{z}_i)\Big)(2t^2)^{m-1}.$$

Therefore, $\mathcal{X}_2(K)$ is defined by
\begin{align*}
z_i+2=\lambda+\lambda^{-1}&=\tau\ne 2, t^2-2, \qquad 1\le i\le m,   \\
H(\lambda,t;r_1,\ldots,r_m)&=0,
\end{align*}
which are $m+1$ equations in the $m+2$ variables $\lambda,t,r_1,\ldots,r_m$. Remember that $z_i,\dot{z}_i,\grave{z}_i$ are polynomials in $t,r_i$.

\begin{rmk}
\rm Notice that (\ref{eq:generic-1}) remains true if $\lambda$ is replaced by $\lambda^{-1}$ (and the condition (\ref{eq:lambda}) is discarded).
With some more efforts, actually we may further rewrite (\ref{eq:generic-1}) as a polynomial in $\tau$, $t$, $r_1,\ldots,r_m$.

To keep things as explicit as possible, we prefer to write (\ref{eq:generic-0}) and (\ref{eq:tau=-2}).
\end{rmk}

\subsection{$\tau=t^2-2\ne \pm 2$}\label{sec:X'}

Abbreviate $\vartheta^{\rm ne}_{T_i}$, $\vartheta^{\rm sw}_{T_i}$ to $\vartheta^{\rm ne}_i$, $\vartheta^{\rm sw}_{i}$, respectively.

Let $t=\kappa+\kappa^{-1}$, with $|\kappa|>1$ or $|\kappa|=1$, $\arg\kappa\in(0,\pi/2)\cup(\pi/2,\pi)$. Up to conjugacy we may assume $\mathbf{g}=\mathbf{d}(\kappa^2)$. For each $i$, either $\mathbf{x}_i^{\rm nw}=\mathbf{x}_i^{\rm ne}=\mathbf{d}(\kappa)$, or there exists $\xi_i\ne 0$ such that
\begin{align}
\mathbf{x}_i^{\rm nw}=\mathbf{u}_\kappa^{\varepsilon_i}(\xi_i), \qquad \mathbf{x}_i^{\rm ne}=\mathbf{u}_\kappa^{\varepsilon_i}(-\kappa^{-2\varepsilon_i}\xi_i),  \label{eq:excptional}
\end{align}
with $\varepsilon_i\in\{+,-\}$, and $\pm$ identified with $\pm 1$. Set $\varepsilon_i=0$ when $\mathbf{x}_i^{\rm nw}=\mathbf{x}_i^{\rm ne}=\mathbf{d}(\kappa)$.

For each $i$, there are five possibilities.
\begin{enumerate}
  \item[\rm 1)] If $\varepsilon_i=\varepsilon_{i+1}=0$, then $\rho|_{T_i}$ is abelian.
  \item[\rm 2)] If $\varepsilon_i=0\ne\varepsilon_{i+1}$, then $\rho|_{T_i}$ is nonabelian reducible.
       By Lemma \ref{lem:reducible2}, $\mathbf{x}_i^{\rm ne}=\mathbf{x}_i^{\rm nw}$ is equivalent to
       \begin{align}
       \vartheta_{i}^{\rm ne}(\kappa^{\varepsilon_{i+1}})=0.  \label{eq:X'-Eq1}
       \end{align}
  \item[\rm 3)] If $\varepsilon_i\ne0=\varepsilon_{i+1}$, the situation is similar, but this time the constraint is
       \begin{align}
       \vartheta_{i}^{\rm ne}(\kappa^{\varepsilon_i})=(1+\kappa^{-2\varepsilon_i})\vartheta_{i}^{\rm sw}(\kappa^{\varepsilon_i}).  \label{eq:X'-Eq2}
       \end{align}
  \item[\rm 4)] If $\varepsilon_i\varepsilon_{i+1}=+$, then $\rho|_{T_i}$ is nonabelian reducible, and
        $$\vartheta_i^{\rm ne}(\kappa^{\varepsilon_{i}})\ne 0, \qquad \vartheta_{i}^{\rm ne}(\kappa^{\varepsilon_i})\ne(1+\kappa^{-2\varepsilon_i})\vartheta_{i}^{\rm sw}(\kappa^{\varepsilon_i}).$$

        Suppose $\mathbf{y}_i=\mathbf{u}^{\varepsilon_i}_\kappa(\zeta_i)$. By the proof of Lemma \ref{lem:reducible2},
        \begin{align*}
        \big(1-\vartheta_{i}^{\rm ne}(\kappa^{\varepsilon_i})\big)\xi_i+\vartheta_{i}^{\rm ne}(\kappa^{\varepsilon_i})\zeta_i&=-\kappa^{-2\varepsilon_i}\xi_i, \\
        \big(1-\vartheta_{i}^{\rm sw}(\kappa^{\varepsilon_i})\big)\xi_i+\vartheta_{i}^{\rm sw}(\kappa^{\varepsilon_i})\zeta_i&=\xi_{i+1}.
        \end{align*}
        Eliminating $\zeta_i$ and noting that $\vartheta^{\rm ne}_{i}(\kappa^{\varepsilon_i})\ne 0$, we obtain
        \begin{align}
        \frac{\xi_{i+1}}{\xi_i}=1-(1+\kappa^{-2\varepsilon_i})\frac{\vartheta_{i}^{\rm sw}(\kappa^{\varepsilon_i})}{\vartheta^{\rm ne}_{i}(\kappa^{\varepsilon_i})}. \label{eq:X'-Eq3}
        \end{align}
  \item[\rm 5)] If $\varepsilon_i\varepsilon_{i+1}=-$, then
        \begin{align}
        z_i=t^2-4, \qquad \dot{z}_i=-\xi_i\xi_{i+1}, \qquad \grave{z}_i=\kappa^{2\varepsilon_i}\xi_i\xi_{i+1}, \label{eq:X'-Eq4}
        \end{align}
        Hence $\grave{z}_i=-\kappa^{-2\varepsilon_i}\dot{z}_i\ne 0$.

        In this case, $\rho|_{T_i}$ is irreducible, so that $r_i\ne 2, t^2-2$.

        We emphasize that $r_i$ should be a common root of $z_i=t^2-4$ and $\grave{z}_i=-\kappa^{2\varepsilon_i}\dot{z}_i$, although sometimes $z_i=t^2-4$ follows from $\grave{z}_i=-\kappa^{2\varepsilon_i}\dot{z}_i$.
\end{enumerate}

\begin{rmk} \label{rmk:X'-II}
\rm The irreducibility of $\rho$ is equivalent to $\{+,-\}\subseteq\{\varepsilon_1,\ldots,\varepsilon_m\}$.

If $K$ is even, then $\rho|_{T_m}$ has been guaranteed to be irreducible, so that $\varepsilon_m\varepsilon_{1}=-$. To see this, note that $\mathbf{x}_m^{\rm ne}$ can always be conjugated to $\mathbf{x}_m^{-1}$ via elements in ${\rm Im}(\rho|_{T_m})$; if $\rho|_{T_m}$ is irreducible, then (\ref{eq:excptional}) would not hold for $i=m$.
\end{rmk}

Let
$$\Xi=\big\{\vec{\varepsilon}=(\varepsilon_1,\ldots,\varepsilon_m)\in\{0,+,-\}^m\colon \{+,-\}\subseteq\{\varepsilon_1,\ldots,\varepsilon_m\}\big\}.$$
Associated to each $\vec{\varepsilon}\in\Xi$ is a system of equations in variables $\kappa,r_1,\ldots,r_m$ and the $\xi_i$'s for  $\varepsilon_i\ne 0$, as (\ref{eq:X'-Eq1})--(\ref{eq:X'-Eq4}). Each solution defines a representation of $K$, for which we can fix the conjugacy indeterminacy via a diagonal matrix so as to set exactly one $\xi_i=1$.
Denote the set of the characters of these representations by $\mathcal{X}'(K;\vec{\varepsilon})$.
Then
$$\mathcal{X}'(K)=\bigsqcup_{\vec{\epsilon}\in\Xi}\mathcal{X}'(K;\vec{\epsilon}).$$

Note that whenever $\varepsilon_\ell=0$ for some $\ell$, there must exist some $i,j\in\{1,\ldots,m\}$ and $\iota,\epsilon\in\{\pm1\}$ such that
\begin{align}
\vartheta_i^{\rm ne}(\kappa^{\iota})=0, \qquad \vartheta_j^{\rm ne}(\kappa^{\epsilon})=(1+\kappa^{-2\epsilon})\vartheta_j^{\rm sw}(\kappa^{\epsilon}). \label{eq:non-generic}
\end{align}
This imposes a numerical constraint on $[p_i/q_i]$ and $[p_j/q_j]$.

If $\varepsilon_i\ne 0$ for all $i$, then
\begin{align*}
1=\prod_{i=1}^m\xi_{i+1}^{\varepsilon_{i+1}}\xi_i^{-\varepsilon_i}=\prod_{i=1}^mf_i^{\varepsilon_i}, 
\end{align*}
with
\begin{align*}
f_i=\begin{cases} 1-(1+\kappa^{-2\varepsilon_i})\vartheta_{i}^{\rm sw}(\kappa^{\varepsilon_i})/\vartheta^{\rm ne}_{i}(\kappa^{\varepsilon_i}), &\varepsilon_i\varepsilon_{i+1}=+, \\
-1/\dot{z}_i, &\varepsilon_i\varepsilon_{i+1}=-. \end{cases}
\end{align*}

Thus, we establish
\begin{prop}\label{prop:X'}
Suppose that for any $i,j$ and any $\iota,\epsilon\in\{\pm1\}$, the equations in (\ref{eq:non-generic}) have no common solution. 
Then $\mathcal{X}'(K)$ consists of finitely many points.
\end{prop}

\vspace{2mm}

Haimiao Chen, Department of Mathematics, Beijing Technology and Business University, 11\# Fucheng Road, Haidian District, Beijing, China.

E-mail: \emph{chenhm@math.pku.edu.cn}

\end{document}